\documentclass[11pt]{article}
\usepackage{amsmath}
\usepackage{amssymb}
\usepackage{amsthm}
\usepackage[usenames]{color}
\usepackage{amscd}
\usepackage[colorlinks=true,linkcolor=webgreen,filecolor=webbrown,
citecolor=webgreen]{hyperref}
\definecolor{webgreen}{rgb}{0,.5,0}
\definecolor{webbrown}{rgb}{.6,0,0}

\def\C{{\mathbb{C}}}

\def\N{{\mathbb{N}}}
\def\Z{{\mathbb{Z}}}
\def\1{{\bf 1}}

\def\lcm{\operatorname{lcm}}

\newtheorem{theorem}{Theorem}
\newtheorem{cor}{Corollary}

\begin{document}

\title{{\bf Representing and counting the subgroups of the group $\Z_m \times \Z_n$}}
\author{Mario Hampejs, Nicki Holighaus, L\'aszl\'o T\'oth, and Christoph Wiesmeyr}
\date{}
\maketitle

\centerline{Journal of Numbers, vol. 2014, Article ID 491428}
\centerline{\url{http://dx.doi.org./10.1155/2014/491428}}

\begin{abstract} We deduce a simple representation and the invariant factor decompositions of the
subgroups of the group $\Z_m \times \Z_n$, where $m$ and $n$ are
arbitrary positive integers. We obtain formulas for the total number of
subgroups and the number of subgroups of a given order.
\end{abstract}

{\sl 2010 Mathematics Subject Classification}:  20K01, 20K27, 05A15, 11A25

{\it Key Words and Phrases}: cyclic group, direct product, finite Abelian
group of rank two, subgroup, number of subgroups, multiplicative arithmetic function

\section{Introduction}

Let $\Z_m$ be the group of residue classes modulo $m$ and consider
the direct product $G=\Z_m \times \Z_n$, where $m$ and $n$ are
arbitrary positive integers. This paper aims to deduce a simple
representation and the invariant factor decompositions of the
subgroups of the group $G$. As consequences we derive formulas for
the number of certain types of subgroups of $G$, including the total
number $s(m,n)$ of its subgroups and the number $s_k(m,n)$ of its
subgroups of order $k$ ($k\mid mn$).

Subgroups of $\Z \times \Z$ (sublattices of the two dimensional
integer lattice) and associated counting functions were considered
by several authors in pure and applied mathematics. It is known, for
example, that the number of subgroups of index $n$ in $\Z \times \Z$
is $\sigma(n)$, the sum of the (positive) divisors of $n$. See,
e.g., \cite{Gra2005}, \cite{Zou2006}, \cite[item A001615]{OEIS}.
Although features of the subgroups of $G$ are not only interesting
by their own but have also applications, one of them described
below, it seems that a synthesis on subgroups of $G$ can not be
found in the literature.

In the case $m=n$ the subgroups of $\Z_n \times \Z_n$ play an
important role in numerical harmonic analysis, more specifically in
the field of applied time-frequency analysis. Time-frequency
analysis attempts to investigate function behavior via a phase space
representation given by the short-time Fourier transform
\cite{Gro2001}. The short-time Fourier coefficients of a function
$f$ are given by inner products with translated modulations (or
time-frequency shifts) of a prototype function $g$, assumed to be
well-localized in phase space, e.g., a Gaussian. In applications, the
phase space corresponding to discrete, finite functions (or vectors)
belonging to $\C^n$ is exactly $\Z_n \times \Z_n$. Concerned with
the question of reconstruction from samples of short-time Fourier
transforms, it has been found that when sampling on lattices, i.e.,
subgroups of $\Z_n \times \Z_n$, the associated analysis and
reconstruction operators are particularly rich in structure, which,
in turn, can be exploited for efficient implementation,
cf.~\cite{KutStr2005}, \cite{Lee2001}, \cite{Str1998} and references
therein. It is of particular interest to find subgroups in a certain
range of cardinality, therefore a complete characterization of these
groups helps choosing the best one for the desired application.

We recall that a finite Abelian group of order $>1$ has rank $r$ if it is
isomorphic to $\Z_{n_1} \times \cdots \times \Z_{n_r}$, where
$n_1,\ldots,n_r \in \N\setminus \{1\}$ and $n_j \mid n_{j+1}$ ($1\le
j\le r-1$), which is the invariant factor decomposition of the given group.
Here the number $r$ is uniquely determined and represents
the minimal number of generators of the group. For
general accounts on finite Abelian groups see, for example, \cite{Mac2012}, \cite{Rot1995}.

It is known that for every finite Abelian group the problem of
counting all subgroups and the subgroups of a given order reduces to
$p$-groups, which follows from the properties of the subgroup
lattice of the group (see \cite{Sch1994}, \cite{Suz1951}). In
particular, for $G=\Z_m \times \Z_n$ this can be formulated as
follows. Assume that $\gcd(m,n)>1$. Then $G$ is an Abelian group of
rank two, since $G\simeq \Z_u \times \Z_v$, where $u=\gcd(m,n)$,
$v=\lcm(m,n)$. Let $u=p_1^{a_1}\cdots p_r^{a_r}$ and
$v=p_1^{b_1}\cdots p_r^{b_r}$ be the prime power factorizations of
$u$ and $v$, respectively, where $0\le a_j\le b_j$ ($1\le j\le r$).
Then
\begin{equation} \label{s_decomp}
s(m,n)=\prod_{j=1}^r s(p_j^{a_j},p_j^{b_j}),
\end{equation}
and
\begin{equation} \label{s_decomp_ord}
s_k(m,n)=\prod_{j=1}^r s_{k_j}(p_j^{a_j},p_j^{b_j}),
\end{equation}
where $k=k_1\cdots k_r$ and $k_j=p_j^{c_j}$ with some exponents
$0\le c_j\le a_j+b_j$ ($1\le j\le r$).

Now consider the $p$-group $\Z_{p^a}\times \Z_{p^b}$, where $0\le
a\le b$. This is of rank two for $1\le a\le b$. One has the simple
explicit formulae:
\begin{equation} \label{prime_pow}
s(p^a,p^b)=
\frac{(b-a+1)p^{a+2}-(b-a-1)p^{a+1}-(a+b+3)p+(a+b+1)}{(p-1)^2},
\end{equation}
\begin{equation} \label{prime_pow_ord}
s_{p^c}(p^a,p^b)= \begin{cases} \frac{p^{c+1}-1}{p-1}, & c\le a\le
b,\\ \frac{p^{a+1}-1}{p-1}, & a\le c\le b,
\\ \frac{p^{a+b-c+1}-1}{p-1},
& a\le b\le c\le a+b.
\end{cases}
\end{equation}

Formula \eqref{prime_pow} was derived by G.~C\u{a}lug\u{a}reanu
\cite[Sect.\ 4]{Cal2004} and recently by  J.~Petrillo \cite[Prop.\
2]{Pet2011} using Goursat's lemma for groups. M.~T\u{a}rn\u{a}uceanu
\cite[Prop.\ 2.9]{Tar2007}, \cite[Th.\ 3.3]{Tar2010} deduced
\eqref{prime_pow} and \eqref{prime_pow_ord} by a method based on
properties of certain attached matrices.

Therefore, $s(m,n)$ and $s_k(m,n)$ can be computed using
\eqref{s_decomp}, \eqref{prime_pow} and \eqref{s_decomp_ord},
\eqref{prime_pow_ord}, respectively. We deduce other formulas for
$s(m,n)$ and $s_k(m,n)$ (Theorems \ref{Th_number_subgroups} and
\ref{Th_number_subgroups_order}), which generalize \eqref{prime_pow}
and \eqref{prime_pow_ord}, and put them in more compact forms. These
are consequences of a simple representation of the subgroups of
$G=\Z_m \times \Z_n$, given in Theorem \ref{Th_repr}. This representation might be known, but
the only source we could find is the paper \cite{KutStr2005}, where only a special case is treated in a different form.
More exactly, in \cite[Lemma 4.1]{KutStr2005} a representation for lattices in $\Z_n \times \Z_n$ of redundancy $2$,
that is subgroups of $\Z_n \times \Z_n$ having index $n/2$ is given, using matrices in Hermite normal form.
Theorem \ref{Th_decomp} gives the invariant factor
decompositions of the subgroups of $G$. We also consider the number
of cyclic subgroups of $\Z_m \times \Z_n$ (Theorem
\ref{Th_number_cyclic_subgroups}) and the number of subgroups of a
given exponent in $\Z_n \times \Z_n$ (Theorem
\ref{Th_number_subgroups_exponent}).

Our approach is elementary, using only simple group-theoretic and number-theoretic arguments.
The proofs are given in Section \ref{Sect_4}.

Throughout the paper we use the notations: $\N=\{1,2,\ldots\}$, $\N_0=\{0,1,2,\ldots\}$,
$\tau(n)$ and $\sigma(n)$ are the number and the sum, respectively, of the
positive divisors of $n$, $\psi(n)=n\prod_{p\mid n} (1+1/p)$ is the
Dedekind function, $\omega(n)$ stands for the number of distinct prime factors of $n$,
$\mu$ is the M\"{o}bius function, $\phi$ denotes Euler's totient function, $\zeta$ is the
Riemann zeta function.


\section{Subgroups of $\Z_m \times \Z_n$} \label{Sect_2}

The subgroups of $\Z_m \times \Z_n$ can be identified and visualized
in the plane with sublattices of the lattice $\Z_m \times \Z_n$.
Every two dimensional sublattice is generated by two basis vectors.
For example, the Figure shows the subgroup of $\Z_{12} \times
\Z_{12}$ having the basis vectors $(3,0)$ and $(1,2)$.

\medskip
\centerline{$
\begin{array}{ccccccccccccc}
11 & . & . & . & . & . & . & . & . & . & . & . & . \\
10 & . & . & \bullet & . & . & \bullet & . & . & \bullet & . & . &
\bullet
\\
9 & . & . & . & . & . & . & . & . & . & . & . & . \\
8 & . & \bullet & . & . & \bullet & . & . & \bullet & . & . & \bullet & . \\
7 & . & . & . & . & . & . & . & . & . & . & . & . \\
6 & \bullet & . & . & \bullet & . & . & \bullet & . & . & \bullet & . & . \\
5 & . & . & . & . & . & . & . & . & . & . & . & . \\
4 & . & . & \bullet & . & . & \bullet & . & . & \bullet & . & . & \bullet \\
3 & . & . & . & . & . & . & . & . & . & . & . & . \\
2 & ._{b} & \bullet & . & . & \bullet & . & . & \bullet & . & . &
\bullet & . \\
1 & . & . & . & . & . & . & . & . & . & . & . & . \\
0 & \bullet & ._{s} & . & \bullet _{a} & . & . & \bullet & . & . &
\bullet &
. & . \\
& 0\text{ } & 1\text{ } & 2\text{ } & 3\text{ } & 4\text{ } & 5\text{ } & 6%
\text{ } & 7\text{ } & 8\text{ } & 9\text{ } & 10 & 11%
\end{array}
$} \vskip1mm

\centerline{Figure}

\medskip

This suggests the following representation of the subgroups:

\begin{theorem} \label{Th_repr} For every $m,n\in \N$ let
\begin{equation}
I_{m,n}:=\{(a,b,t)\in \N^2\times \N_0: a\mid m, b\mid n, 0\le t\le \gcd(a,n/b)-1\}
\end{equation}
and for $(a,b,t)\in I_{m,n}$ define
\begin{equation} \label{6_repr}
H_{a,b,t}:= \{(ia+jta/\gcd(a,n/b), jb): 0\le i\le m/a-1, 0\le
j\le n/b-1\}.
\end{equation}

Then $H_{a,b,t}$ is a subgroup of order $\frac{mn}{ab}$ of $\Z_m \times \Z_n$
and the map $(a,b,t)\mapsto H_{a,b,t}$ is a bijection between the
set $I_{m,n}$ and the set of subgroups of $\Z_m \times \Z_n$.
\end{theorem}

Note that for the subgroup $H_{a,b,t}$ the basis vectors mentioned
above are $(a,0)$ and $(s,b)$, where
\begin{equation}
s=\frac{ta}{\gcd(a,n/b)}.
\end{equation}

This notation for $s$ will be used also in the rest of the paper.
Note also that in the case $a\ne m$, $b\ne n$ the area of the
parallelogram spanned by the basis vectors is $ab$, exactly the
index of $H_{a,b,t}$.

We say that a subgroup $H=H_{a,b,t}$ is a subproduct of $\Z_m \times
\Z_n$ if $H=H_1\times H_2$, where $H_1$ and $H_2$ are subgroups of
$\Z_m$ and $\Z_n$, respectively.

\begin{theorem} \label{Th_decomp}

i) The invariant factor decomposition of the subgroup $H_{a,b,t}$ is
given by
\begin{equation}
H_{a,b,t} \simeq  \Z_{\alpha} \times \Z_{\beta},
\end{equation}
where
\begin{equation}
\alpha=\gcd(m/a,n/b,ns/(ab)), \quad \beta =\frac{mn}{ab\alpha}
\end{equation}
satisfying $\alpha \mid \beta$.

ii) The exponent of the subgroup $H_{a,b,t}$ is $\beta$.

iii) The subgroup $H_{a,b,t}$ is cyclic if and only if $\alpha=1$.

iv) The subgroup $H_{a,b,t}$ is a subproduct if and only if $t=0$
and $H_{a,b,0}=\Z_{m/a}\times \Z_{n/b}$. Here $H_{a,b,0}$ is cyclic
if and only if $\gcd(m/a,n/b)=1$.
\end{theorem}

For example, for the subgroup represented by the Figure one has
$m=n=12$, $a=3$, $b=2$, $s=1$, $\alpha=2$, $\beta=12$, and this subgroup is
isomorphic to $\Z_2\times \Z_{12}$. It is not cyclic and is not a
subproduct.

According to Theorem \ref{Th_repr}, the number $s(m,n)$ of subgroups of $\Z_m
\times \Z_n$ can be obtained by counting the elements of the set
$I_{m,n}$. We deduce

\begin{theorem} \label{Th_number_subgroups} For every $m,n\in \N$, $s(m,n)$
is given by
\begin{equation} \label{total_number_subgroups}
s(m,n)= \sum_{a\mid m, b\mid n} \gcd(a,b)
\end{equation}
\begin{equation} \label{total_number_subgroups_var_1}
= \sum_{d\mid \gcd(m,n)} \phi(d)\tau(m/d)\tau(n/d)
\end{equation}
\begin{equation} \label{total_number_subgroups_var_2}
= \sum_{d\mid \gcd(m,n)} d\tau(mn/d^2).
\end{equation}
\end{theorem}

Formula \eqref{total_number_subgroups} is a special case of a formula
representing the number of all subgroups of a
class of groups formed as cyclic extensions of cyclic groups,
deduced by W.~C.~Calhoun \cite{Cal1987} and having a laborious
proof. Note that formula \eqref{total_number_subgroups} is
given, without proof in \cite[item A054584]{OEIS}.

Note also that the function $(m,n)\mapsto s(m,n)$ is representing a
multiplicative arithmetic function of two variables, that is,
$s(mm',nn')= s(m,n)s(m',n')$ holds for any $m,n,m',n'\in \N$ such
that $\gcd(mn,m'n')=1$. This property, which is in concordance with
\eqref{s_decomp}, is a direct consequence of formula
\eqref{total_number_subgroups}. See Section \ref{Sect_5}.

Let $N(a,b,c)$ denote the number of solutions $(x,y,z,t)\in
\N^4$ of the system of equations $xy=a, zt=b, xz=c$.

\begin{theorem} \label{Th_number_subgroups_order} For every $k,m,n\in \N$ such that $k\mid mn$,
\begin{equation} \label{total_number_subgroups_ord_k}
s_k(m,n)= \sum_{\substack{a\mid m, b\mid n\\ mb/a=k}} \gcd(a,b)
\end{equation}
\begin{equation} \label{total_number_subgroups_ord_k_var_1}
=\sum_{\substack{d\mid \gcd(k,m) \\ e\mid \gcd(k,n)\\ k\mid de}} \phi(de/k)
\end{equation}
\begin{equation} \label{total_number_subgroups_ord_k_var_2}
= \sum_{d\mid \gcd(m,n,k)} \phi(d) N(m/d, n/d, k/d).
\end{equation}
\end{theorem}

The identities \eqref{prime_pow} and \eqref{prime_pow_ord} can be
easily deduced from each of the identities given in Theorems \ref
{Th_number_subgroups} and \ref{Th_number_subgroups_order},
respectively.

\begin{theorem} \label{Th_number_cyclic_subgroups} Let $m,n\in \N$.

i) The number $c(m,n)$ of cyclic subgroups of $\Z_m \times \Z_n$ is given by
\begin{equation} \label{total_number_cyclic_subgroups}
c(m,n)= \sum_{\substack{a\mid m, b\mid n\\ \gcd(m/a,n/b)=1}} \gcd(a,b)
\end{equation}
\begin{equation} \label{total_number_cyclic_subgroups_var_1}
= \sum_{a\mid m, b\mid n} \phi(\gcd(a,b))
\end{equation}
\begin{equation} \label{total_number_cyclic_subgroups_var_2}
= \sum_{d\mid \gcd(m,n)} (\mu*\phi)(d)\tau(m/d)\tau(n/d)
\end{equation}
\begin{equation} \label{total_number_cyclic_subgroups_var_3}
= \sum_{d\mid \gcd(m,n)} \phi(d) \tau(mn/d^2).
\end{equation}

ii) The number of subproducts of $\Z_m \times \Z_n$ is
$\tau(m)\tau(n)$ and the number of its  cyclic subproducts is
$\tau(mn)$.
\end{theorem}

Formula \eqref{total_number_cyclic_subgroups_var_1}, as a special
case of an identity valid for arbitrary finite Abelian groups, was
derived by the third author \cite{Tot2011,Tot2012} using different
arguments. The function $(m,n)\mapsto c(m,n)$ is also
multiplicative.


\section{Subgroups of $\Z_n \times \Z_n$} \label{Sect_3}

In the case $m=n$, which is of special interest in applications, the
results given in the previous section can be easily used. We point
out that $n\mapsto s(n):=s(n,n)$ and $n\mapsto c(n):=c(n,n)$ are
multiplicative arithmetic functions of a single variable (sequences
\cite[items A060724, A060648]{OEIS}). They can be written in the
form of Dirichlet convolutions as shown by the next Corollaries.

\begin{cor} \label{Cor_1} For every $n\in \N$,
\begin{equation}
s(n)=  \sum_{de=n} \phi(d) \tau^2(e)
\end{equation}
\begin{equation} \label{repr_var}
= \sum_{de=n} d\, \tau(e^2).
\end{equation}
\end{cor}

\begin{cor} \label{Cor_1} For every $n\in \N$,
\begin{equation} \label{omega_formula}
c(n)=  \sum_{de=n} d\, 2^{\omega(e)}
\end{equation}
\begin{equation} \label{repr_var}
= \sum_{de=n} \phi(d) \tau(e^2).
\end{equation}
\end{cor}

Further convolutional representations can also be given, for example,
\begin{equation} \label{connect_s_c}
s(n)= \sum_{de=n} \tau(d)\psi(e), \qquad  c(n)=  \sum_{d\mid n} \psi(d),
\end{equation}
all of these follow from the Dirichlet-series representations
\begin{equation} \label{Dir_ser_s}
\sum_{n=1}^{\infty} \frac{s(n)}{n^z}=
\frac{\zeta^3(z)\zeta(z-1)}{\zeta(2z)},
\end{equation}
\begin{equation} \label{Dir_ser_c}
\sum_{n=1}^{\infty} \frac{c(n)}{n^z}=
\frac{\zeta^2(z)\zeta(z-1)}{\zeta(2z)},
\end{equation}
valid for $z\in \C$, $\Re (z) >2$.

Observe that
\begin{equation*}
s(n)= \sum_{d\mid n} c(d) \quad (n\in \N),
\end{equation*}
which is a simple consequence of \eqref{connect_s_c} or of
\eqref{Dir_ser_s} and \eqref{Dir_ser_c}. It also follows from the
next result.

\begin{theorem} \label{Th_number_subgroups_exponent}
For every $n,\delta\in \N$ with $\delta \mid n$ the number of
subgroups of exponent  $\delta$ of $\Z_n \times \Z_n$
equals the number of cyclic subgroups of $\Z_{\delta} \times \Z_{\delta}$.
\end{theorem}


\section{Proofs} \label{Sect_4}

\begin{proof} (for Theorem \ref{Th_repr}) Let $H$ be a subgroup of $G=\Z_m \times \Z_n$.
Consider the natural projection
$\pi_2:G\to \Z_n$ given by $\pi_2(x,y)=y$. Then $\pi_2(H)$ is a
subgroup of $\Z_n$ and there is a unique divisor $b$ of $n$ such
that $\pi_2(H)=\langle b \rangle:=\{jb: 0\le j\le n/b-1\}$. Let
$s\ge 0$ be minimal such that $(s,b)\in H$.

Furthermore, consider the natural inclusion $\iota_1:\Z_m\to G$
given by $\iota_1(x)=(x,0)$. Then $\iota_1^{-1}(H)$ is a subgroup of
$\Z_m$ and there exists a unique divisor $a$ of $m$ such that
$\iota_1^{-1}(H)=\langle a \rangle$.

We show that $H=\{(ia+js,jb): i,j\in \Z \}$. Indeed, for every
$i,j\in \Z$, $(ia+js,jb)= i(a,0)+j(s,b)\in H$. On the other hand,
for every $(u,v)\in H$ one has $v\in \pi_2(H)$ and hence there is
$j\in \Z$ such that $v=jb$. We obtain $(u-js,0)= (u,v)- j(s,b)\in
H$, $u-js\in \iota_1^{-1}(H)$ and there is $i\in \Z$ with
$u-js=ia$.

Here a necessary condition is that $(sn/b,0)\in H$
(obtained for $i=0$, $j=n/b$), that is $a\mid sn/b$, equivalent to
$a/\gcd(a,n/b) \mid s$. Clearly, if this is verified, then for the above representation of $H$ it is
enough to take the values $0\le i\le m/a-1$ and $0\le j\le n/b-1$.

Also, dividing $s$ by $a$ we have $s=aq+r$ with $0\le r<a$ and $(r,b)=(s,b)-q(a,0)\in H$,
showing that $s<a$, by its minimality. Hence
$s=ta/\gcd(a,n/b)$ with $0\le t\le \gcd(a,n/b)-1$. Thus we obtain the given representation.

Conversely, every $(a,b,t)\in I_{m,n}$ generates a subgroup $H_{a,b,t}$
of order $mn/(ab)$ of $\Z_m \times \Z_n$ and the proof is complete.
\end{proof}

\begin{proof} (for Theorem \ref{Th_decomp}) i)-ii) We first determine the exponent of the subgroup
$H_{a,b,t}$. $H_{a,b,t}$ is generated by $(a,0)$ and $(s,b)$, hence its exponent is the
least common multiple of the orders of these two elements. The order
of $(a,0)$ is $m/a$. To compute the order of $(s,b)$ note that
$m\mid rs$ if and only if $m/\gcd(m,s) \mid r$. Thus the order of
$(s,b)$ is $\lcm(m/\gcd(m,s),n/b)$. We deduce that the exponent of
$H_{a,b,t}$ is
\begin{equation*}
\lcm \left(\frac{m}{a}, \frac{m}{\gcd(m,s)},
\frac{n}{b}\right) = \lcm \left(\frac{mn}{na},
\frac{mn}{n\gcd(m,s)}, \frac{mn}{mb} \right)
\end{equation*}
\begin{equation*}
=\frac{mn}{\gcd(na,nm,ns,mb)}= \frac{mn}{\gcd(mb,na,ns)} =\beta
\end{equation*}

For every finite Abelian group the rank of a nontrivial subgroup is
at most the rank of the group. Therefore, the rank of $H_{a,b,t}$ is
$1$ or $2$. That is, $H_{a,b,t}\simeq \Z_A \times \Z_B$ with certain
$A,B\in \N$ such that $A \mid B$. Here the exponent of $H_{a,b,t}$ equals that of $\Z_A \times \Z_B$, which is
$\lcm(A,B)=B$. Using ii) already proved we deduce that $B=\beta$.
Since the order of $H_{a,b,t}$ is $AB=mn/(ab)$ we have
$A=mn/(ab\beta)=\alpha$.

iii) According to i), $H_{a,b,t}\simeq \Z_{\alpha} \times \Z_{\beta}$, where $\alpha \mid
\beta$. Hence $H_{a,b,t}$ is cyclic if and only if $\alpha=1$.

iv) The subgroups of $\Z_m$ are of form $\{ia: 0\le i\le m/a-1\}$, where $a\mid m$, and the properties follow from \eqref{6_repr} and iii).
\end{proof}

\begin{proof} (for Theorem \ref{Th_number_subgroups}) By its definition, the number of
elements of the set $I_{m,n}$ is
\begin{equation*}
\sum_{a\mid m, b\mid n} \sum_{0\le t\le \gcd(a,n/b)-1} 1 =\sum_{a\mid
m, b\mid n} \gcd(a,n/b) = \sum_{a\mid m, b\mid n} \gcd(a,b),
\end{equation*}
representing $s(m,n)$. This is formula
\eqref{total_number_subgroups}.

To obtain formula \eqref{total_number_subgroups_var_1} apply the
Gauss formula $n=\sum_{d\mid n} \phi(d)$ ($n\in \N$) by writing:
\begin{equation*}
s(m,n)= \sum_{a\mid m, b\mid n} \sum_{d\mid \gcd(a,b)} \phi(d)=
\sum_{\substack{ax=m\\ by=n}} \sum_{\substack{di=a\\
dj=b}} \phi(d) = \sum_{\substack{dix=m\\ djy=n}} \phi(d)
\end{equation*}
\begin{equation*}
=  \sum_{\substack{du=m\\ dv=n}} \phi(d) \sum_{\substack{ix=u \\ jy=v}}
1= \sum_{\substack{du=m\\ dv=n}} \phi(d) \tau(u)\tau(v)
\end{equation*}
\begin{equation*}
= \sum_{d\mid \gcd(m,n)} \phi(d) \tau(m/d)\tau(n/d).
\end{equation*}

Now \eqref{total_number_subgroups_var_2} follows from
\eqref{total_number_subgroups_var_1} by the Busche-Ramanujan
identity (cf. \cite[Ch.\ 1]{McC1986})
\begin{equation*}
\tau(m)\tau(n) = \sum_{d \mid \gcd(m,n)} \tau(mn/d^2) \quad (m,n\in \N).
\end{equation*}
\end{proof}

\begin{proof} (for Theorem \ref{Th_number_subgroups_order}) According to Theorem \ref{Th_repr},
\begin{equation*}
s_k(m,n)= \sum_{\substack{a\mid m, b\mid n\\ mn/ab=k}} \gcd(a,n/b),
\end{equation*}
giving \eqref{total_number_subgroups_ord_k}, which can be written,
by Gauss' formula again, as
\begin{equation} \label{s_k_proof_line}
s_k(m,n)= \sum_{\substack{a\mid m, b\mid n\\ mb/a=k}}
\sum_{\substack{c\mid a, c\mid b}} \phi(c)= \sum_{\substack{cix=m\\
cjy=n \\ cjx=k}} \phi(c)
\end{equation}
\begin{equation*}
=\sum_{\substack{di=m\\ ey=n}} \sum_{\substack{cx=d\\
cj=e\\ cjx=k}} \phi(c),
\end{equation*}
where in the inner sum one has $c=de/k$ and obtain
\eqref{total_number_subgroups_ord_k_var_1}. Now, to get
\eqref{total_number_subgroups_ord_k_var_2} write \eqref{s_k_proof_line} as
\begin{equation*}
s_k(m,n)= \sum_{\substack{cu=m\\ cv=n\\ cw=k}} \phi(c)
\sum_{\substack{ix=u\\ jy=v\\ jx=w}} 1 = \sum_{\substack{cu=m\\
cv=n\\ cw=k}} \phi(c) N(u,v,w)
\end{equation*}
\begin{equation*}\sum_{c\mid
\gcd(m,n,k)} \phi(c) N(m/c,n/c,k/c),
\end{equation*}
and the proof is complete.
\end{proof}

\begin{proof} (for Theorem \ref{Th_number_cyclic_subgroups}) i) According to Theorems \ref{Th_repr} and
\ref{Th_decomp}/iii) and using that $\sum_{d\mid n} \mu(d)=1$ or $0$,
according to $n=1$ or $n>1$,
\begin{equation*}
c(m,n)= \sum_{a\mid m, b\mid n} \sum_{\substack{1\le s\le a\\ ab\mid ns\\ \gcd(m/a,n/b,ns/ab)=1}} 1 =
\sum_{\substack{ax=m \\ by=n}} \sum_{\substack{1\le s\le a\\ ar=ys \\ \gcd(x,y,r)=1}} 1
\end{equation*}
\begin{equation*}
= \sum_{\substack{ax=m\\ by=n}} \sum_{\substack{1\le s\le a\\
ar=ys}} \sum_{e\mid \gcd(x,y,r)} \mu(e) = \sum_{\substack{aei=m\\
bej=n}} \mu(e) \sum_{\substack{1\le s\le a\\ a/\gcd(a,j)\mid s}} 1,
\end{equation*}
where the inner sum is $\gcd(a,j)$. Hence
\begin{equation} \label{formula_cyclic}
c(m,n) = \sum_{\substack{aei=m\\ bej=n}} \mu(e) \gcd(a,j).
\end{equation}

Now regrouping the terms according to $ei=z$ and $be=t$ we obtain
\begin{equation*}
c(m,n) = \sum_{\substack{az=m\\ jt=n}} \gcd(a,j) \sum_{\substack{ei=z\\ be=t}}
\mu(e) = \sum_{\substack{az=m\\ jt=n}}
\gcd(a,j) \sum_{e\mid \gcd(z,t)} \mu(e)
\end{equation*}
\begin{equation*}
= \sum_{\substack{az=m\\ jt=n\\ \gcd(z,t)=1}} \gcd(a,j),
\end{equation*}
which is \eqref{total_number_cyclic_subgroups}.

The next results follow applying Gauss' formula and the
Busche-Ramanujan formula, similar to the proof of Theorem
\ref{Th_number_subgroups}.

ii) For the subproducts $H_{a,b,0}$ the values $a\mid m$ and $b\mid
n$ can be chosen arbitrary and it follows at once that the number of
subproducts is $\tau(m)\tau(n)$. The number of cyclic subproducts is
\begin{equation*}
\sum_{\substack{a\mid m\\ b\mid n\\ \gcd(m/a,n/b)=1}} 1 = \sum_{\substack{ax= m\\ by=n \\ \gcd(x,y)=1}} 1 =
\sum_{\substack{ax= m\\ by=n}} \sum_{e\mid \gcd(x,y)} \mu(e) =
\end{equation*}
\begin{equation*}
=\sum_{\substack{eA=m\\ eB=n}} \mu(e) \tau(A)\tau(B)= \sum_{e\mid
\gcd(m,n)} \mu(e) \tau(m/e)\tau(n/e)= \tau(mn),
\end{equation*}
by the inverse Busche-Ramanujan identity.
\end{proof}

\begin{proof} (for Theorem \ref{Th_number_subgroups_exponent}) According to Theorem \ref{Th_decomp}/ii),
the number of subgroups of exponent $\delta$ of $\Z_n \times \Z_n$ is
\begin{equation*}
E_{\delta}(n)=\sum_{a\mid n, b\mid n} \sum_{\substack{1\le s\le a\\
ab\mid ns\\ n/\gcd(a,b,s)=\delta}} 1 = \sum_{\substack{ax=n \\
by=n}} \sum_{\substack{1\le s\le a\\ ar=ys \\ \gcd(a,b,s)=n/\delta}}
1.
\end{equation*}

Write $a=a_1n/\delta$, $b=b_1n/\delta$, $s=s_1n/\delta$ with $\gcd(a_1,b_1,s_1)=1$. We deduce,
similar to the proof of Theorem \ref{Th_number_cyclic_subgroups}/i) that
\begin{equation*}
E_{\delta}(n)= \sum_{\substack{eix=\delta\\ ejy=\delta}} \mu(e)
\gcd(i,y),
\end{equation*}
which is exactly $c(\delta,\delta)=c(\delta)$, cf.
\eqref{formula_cyclic}.
\end{proof}

\section{Further remarks} \label{Sect_5}

1) As mentioned in the Section \ref{Sect_2} the functions
$(m,n)\mapsto s(m,n)$ and $(m,n)\mapsto c(m,n)$ are multiplicative
functions of two variables. This follows easily from formulae
\eqref{total_number_subgroups} and
\eqref{total_number_cyclic_subgroups_var_1}, respectively. Namely,
according to those formulae $s(m,n)$ and $c(m,n)$ are two variables
Dirichlet convolutions of the functions $(m,n)\mapsto \gcd(m,n)$ and
$(m,n)\mapsto \phi(\gcd(m,n))$, respectively with the constant $1$
function, all multiplicative. Since convolution preserves the
multiplicativity we deduce that $s(m,n)$ and $c(m,n)$ are also
multiplicative. See \cite[Sect.\ 2]{Tot2011} for details.

2) Asymptotic formulas with sharp error terms for the sums
$\sum_{m,n\le x} s(m,n)$ and $\sum_{m,n\le x} c(m,n)$ were given in
the paper \cite{NowTot2013}.

3) For any finite groups $A$ and $B$ a subgroup $C$ of $A\times B$
is cyclic if and only if $\iota_1^{-1}(C)$ and $\iota_2^{-1}(C)$
have coprime orders, where $\iota_1$ and $\iota_2$ are the natural
inclusions (\cite[Th.\ 4.2]{BauSenZve2011}). In the case $A=\Z_m$,
$B=\Z_n$ and $C=H_{a,b,t}$ one has $\# \iota_1^{-1}(C)=m/a$ and $\#
\iota_2^{-1}(C)= \gcd(n/b,ns/ab)$ and the characterization of the
cyclic subgroups $H_{a,b,t}$ given in Theorem \ref{Th_decomp}/iii)
can be obtained also in this way. It turns out that regarding the
sublattice, $H_{a,b,t}$ is cyclic if and only if the numbers of
points on the horizontal and vertical axes, respectively, are
relatively prime. Note that in the case $m=n$ the above condition
reads $n\gcd(a,b,s)=ab$. Thus it is necessary that $n\mid ab$. The
subgroup on the Figure is not cyclic.

4) Note also the next formula for the number of cyclic subgroups of $\Z_n \times \Z_n$,
derived in \cite[Ex.\ 2]{PakWar2009}:
\begin{equation} \label{2_other_1}
c(n)= \sum_{\lcm(d,e)=n} \gcd(d,e) \quad (n\in \N),
\end{equation}
where the sum is over all ordered pairs $(d,e)$ such that
$\lcm(d,e)=n$. For a short direct proof of \eqref{2_other_1} write
$d=\ell a$, $e=\ell b$ with $\gcd(a,b)=1$. Then $\gcd(d,e)=\ell$,
$\lcm(d,e)=\ell ab$ and obtain
\begin{equation*}
\sum_{\lcm(d,e)=n} \gcd(d,e) = \sum_{\substack{\ell ab=n \\
\gcd(a,b)=1}} \ell = \sum_{\ell k=n} \ell \sum_{\substack{ab = k \\
\gcd(a,b)=1}} 1 = \sum_{\ell k=n} \ell \, 2^{\omega(k)} = c(n),
\end{equation*}
according to \eqref{omega_formula}.

5) Every subgroup $K$ of $\Z \times \Z$ has the representation $K=\{
(ia+js,jb): i,j\in \Z\}$, where $0\le s\le a$, $0\le b$ are unique
integers. This follows like in the proof of Theorem \ref{Th_repr}.
Furthermore, in the case $a,b\ge 1$, $0\le s\le a-1$ the index of $K$ is $ab$
and one obtains at once that the number of subgroups $K$ having
index $n$ ($n\in \N$) is $\sum_{ab=n} \sum_{0\le s \le a-1} 1 =
\sum_{ab=n} a= \sigma(n)$, mentioned in the Introduction.

\section{Acknowledgement} N.~Holighaus was partially supported by
the Austrian Science Fund (FWF) START-project FLAME (Y551-N13).
L.~T\'oth gratefully acknowledges support from the Austrian Science
Fund (FWF) under the project Nr. M1376-N18. C.~Wiesmeyr was
partially supported by EU FET Open grant UNLocX (255931).


\vskip4mm

\noindent M.~Hampejs \\ NuHAG, Faculty of Mathematics, University of Vienna, Oskar Morgenstern
Platz 1, A-1090 Vienna, Austria \\ E-mail:
mario.hampejs@univie.ac.at

\medskip

\noindent N.~Holighaus \\ Acoustics Research Institute, Austrian
Academy of Sciences \\ Wohllebengasse 12-14, A-1040 Vienna, Austria
\\ E-mail: nicki.holighaus@univie.ac.at

\medskip

\noindent L. T\'oth \\
Department of Mathematics, University of P\'ecs \\ Ifj\'us\'ag u. 6, H-7624 P\'ecs, Hungary
\\ and \\ Institute of Mathematics, University of Natural Resources and Life
Sciences, Gregor Mendel Stra{\ss}e 33, A-1180 Vienna, Austria
\\ E-mail: ltoth@gamma.ttk.pte.hu

\medskip

\noindent C.~Wiesmeyr \\ NuHAG, Faculty of Mathematics, University of Vienna, Oskar Morgenstern
Platz 1, A-1090 Vienna, Austria  \\ E-mail:
christoph.wiesmeyr@univie.ac.at


\begin{thebibliography}{99}

\bibitem{BauSenZve2011} K.~Bauer, D.~Sen, P.~Zvengrowski, A generalized Goursat lemma,
Preprint, arXiv:11009.0024 [math.GR].

\bibitem{Cal1987} W.~C.~Calhoun, Counting the subgroups of some finite groups,
\textit{Amer. Math. Monthly}, {\bf 94} (1987), 54--59.

\bibitem{Cal2004} G.~C\u{a}lug\u{a}reanu, The total number of subgroups of a finite
abelian group, \textit{Sci. Math. Jpn.} {\bf 60} (2004), 157--167.

\bibitem{Gra2005} M.~J.~Grady, A group theoretic approach to a famous partition formula,
{\it Amer. Math. Monthly}, {\bf 112} (2005), 645--651.

\bibitem{Gro2001} K.~Gr\"ochenig, {\it Foundations of Time-Frequency Analysis},
Applied and Numerical Harmonic Analysis, Birkh{\"a}user Boston, 2001.

\bibitem{KutStr2005} G.~Kutyniok, T.~Strohmer, Wilson bases for general time-frequency lattices,
{\it SIAM J. Math. Anal.}, {\bf 37} (2005), 685--711.

\bibitem{Lee2001} A.~J.~van~Leest, {\it Non-separable Gabor schemes. Their Design and Implementation}, PhD thesis,
Tech. Univ. Eindhoven, 2001.

\bibitem{NowTot2013} W.~G.~Nowak, L.~T\'oth, On the average number of subgroups of the
group $\Z_m \times \Z_n$, {\it Int. J. Number Theory}, {\bf 10} (2014), 363-374.

\bibitem{Mac2012} A.~Mach\`{i}, {\it Groups. An Introduction to Ideas and Methods of the Theory of Groups}, Springer, 2012.

\bibitem{McC1986} P.~J.~McCarthy, {\it Introduction to Arithmetical Functions},
Springer, 1986.

\bibitem{PakWar2009} A.~Pakapongpun, T.~Ward, Functorial orbit counting, {\it J. Integer Sequences},
{\bf 12} (2009), Article 09.2.4, 20 pp.

\bibitem{Pet2011} J.~Petrillo, Counting subgroups in a direct product of finite cyclic grups,
{\it College Math J.}, {\bf 42} (2011), 215--222.

\bibitem{Rot1995} J.~J.~Rotman, {\it An Introduction to the Theory of Groups}, Fourth Ed., Springer, 1995.

\bibitem{Sch1994} R.~Schmidt, {\it Subgroup Lattices of Groups}, de Gruyter Expositions
in Mathematics 14, de Gruyter, Berlin, 1994.

\bibitem{Str1998} T.~Strohmer, {\it Numerical algorithms for discrete Gabor expansions}, In H.~G.~
Feichtinger and T.~Strohmer, editors, Gabor Analysis and Algorithms: Theory and Applications, pp. 267--294,
Birkh{\"a}user Boston, 1998.

\bibitem{Suz1951} M.~Suzuki, On the lattice of subgroups of finite
groups, {\it Trans. Amer. Math. Soc.}, {\bf 70} (1951), 345--371.

\bibitem{Tar2007} M.~T\u{a}rn\u{a}uceanu, A new method of proving some classical theorems of abelian
groups, {\it Southeast Asian Bull. Math.}, {\bf 31} (2007),
1191--1203.

\bibitem{Tar2010} M.~T\u{a}rn\u{a}uceanu, An arithmetic method of counting the
subgroups of a finite abelian group, {\it Bull. Math. Soc. Sci.
Math. Roumanie (N.S.)}, {\bf 53(101)} (2010), 373--386.

\bibitem{Tot2011} L.~T\'oth, Menon's identity and arithmetical sums
representing functions of several variables, {\it Rend. Sem. Mat.
Univ. Politec. Torino}, {\bf 69} (2011), 97--110.

\bibitem{Tot2012} L.~T\'oth, On the number of cyclic subgroups of a finite Abelian
group, {\it Bull. Math. Soc. Sci. Math. Roumanie (N.S.)}, {\bf
55(103)} (2012), 423--428.

\bibitem{Zou2006} Y.~M.~Zou, Gaussian binomials and the number of sublattices,
{\it Acta Cryst.}, {\bf 62} (2006), 409--410.

\bibitem{OEIS} The On-Line Encyclopedia of Integer Sequences, {http://oeis.org}.

\end{thebibliography}
\end{document}